\documentclass[12pt]{amsart}
\usepackage{amsaddr}
\usepackage{xcolor}
\usepackage[
top=.7in, bottom=.7in,
left=1in, right=1in]{geometry}
\title{Generalized Carlos Scales}

\author{Andrew V. Sills}
\address{Department of Mathematical Sciences\\ Georgia Southern University\\
Statesboro, Georgia, 30460, USA}
\email{asills@georgiasouthern.edu}
\date{\today}

\theoremstyle{remark}


\begin{document}
\maketitle

\begin{abstract}
In 1986, composer Wendy Carlos introduced three unusual musical scales she called alpha, beta and gamma, equal temperament-inspired scales that de-emphasize the octave as the primary interval in favor of the major and minor thirds and the perfect fifth. 
A derivation of the alpha, beta, and gamma scales due to David Benson is
generalized to produce many Carlos-type scales.  \end{abstract}

\section{Introduction: equal temperament}

For any musical scale where the spaces between available notes are equal, there is a 
positive real number $x_0$, so that multiplying the frequency of a given note by $2^{x_0}$ results in the 
frequency of the next higher available note.   
In such a scale, we will refer to the interval between adjacent pitches
as a \emph{unit}\footnote{The author prefers the term ``unit'' to ``step'' 
since the latter could be confused
with the whole step/half step terminology common in Western music theory.  Our
``unit'' is a generalization of the semi-tone or ``half-step''.}.
 For example, in twelve tone equal temperament (12-ET),
$x_0 = 1/12$ and the unit is the equal-tempered semitone.   Once one has decided that
the octave must be a perfect $2:1$ ratio and that there will be $12$ units per octave, then
$x_0 = 1/12$ is the unique solution.  With this solution, the 12-ET perfect fifth has a ratio of $2^{7/12}$,
which approximates the just perfect fifth of $3:2$ quite well, falling just $1.96$ cents narrower
than just.
The minor third, at $2^{3/12}$, is a less satisfactory approximation of the just ratio $6:5$, falling 
$15.64$ cents narrow.  The major third, at $2^{4/12}$, is similarly not such a great approximation to the 
just ratio of $5:4$, as it is $13.69$ cents wider.

A well trained human ear is perhaps sensitive to deviations in pitch at the $5$ cent level.
Given the centrality of the major and minor triad in Western music, the 12-ET solution may seem
unacceptable to those of us with a fine-tuned sense of pitch. 

\section{Deriving the original Carlos scales and some new Carlos-type scales}
In the 1980s, Wendy Carlos effectively asked the following question: 
what if we design an equal unit scale giving
 priority to the major and minor third and the perfect fifth, so the major and minor triad built on 
the tonic approximates just intonation as closely as possible, and let the octave fall where it
may?  In~\cite{C90}, Carlos reveals that she ``wrote a computer program to perform a
precise deep-search investigation'' that produced three new interesting scales that she
named alpha, beta, and gamma.   These new musical scales were introduced to the
public in Carlos' 1986 album \emph{Beauty in the Beast}~\cite{C86}.
The methodology of the computer search that Carlos used to derive the 
scales is not described in~\cite{C90}, but guided by a method 
given by Benson~\cite[pp. 231--233]{B06},
we can not only re-derive the Carlos alpha, beta, and gamma scales, but also many other
scales guided by the same principle.

We note that perhaps the simplest way of obtaining
the alpha, beta, and gamma scales is by dividing the just perfect fifth (3:2 ratio) into 
equal parts: 9 equal parts for alpha, 11 equal parts for beta, and 20 equal parts for gamma.
This implies the alpha unit is 77.995 cents, the beta unit is 63.814 cents, and the gamma unit is 35.09775 cents.   These values differ ever so slightly from the alpha, beta, and gamma units as
derived by Benson's method (77.965, 63.833, and 35.099 cents respectively). To the human ear, these discrepancies would be imperceptible, but is nonetheless of potential mathematical interest.
By dividing the just perfect fifth, the Carlos scales would, by definition, have a perfectly in tune
fifth, while the Benson method yields scales with ever so slightly wide or narrow fifths.
Carlos~\cite[p. 2]{C90} hints that she intended the perfect fifth to be exactly just, i.e.
an exact $3:2$ ratio, by writing,
``If you try to play through a one octave scale of Alpha, you'd find there are 4 steps to the minor
third, 5 steps to the major third and 9 steps to the perfect (no kidding) fifth, \dots''
Elsewhere in Carlos's writings~\cite{C87,C90}, she uses the values 78.0, 63.8, and 35.1 cents 
respectively for
unit size in the alpha, beta, and gamma scales: thus to three significant figures, there is
no difference between the equal division of the just perfect fifth and Benson's method.
We proceed using Benson's method.

The goal of a ``Carlos-type'' scale is to choose positive integers $a$ and $b$,  with $a<b$, so that 
a minor third is $a$ units above the tonic, the major third is $b$ units above the tonic, and the
perfect fifth is $a+b$ units above the tonic.  Thus we want to find a 
unit size $x$ that simultaneously approximates

\begin{equation} (2^x)^a \approx 6/5, \qquad (2^x)^b \approx 5/4,
\mbox{ and } (2^x)^{a+b} \approx 3/2. 
\label{apr1}
\end{equation}

The approximations given in line~\eqref{apr1} are equivalent to
\begin{equation} 
 \frac{ 2^{ax}}{6/5} \approx 1, \qquad \frac{2^{bx}}{4/3} \approx 1, \qquad \frac{2^{(a+b)x}}{3/2}
 \approx 1,
\label{apr2}
\end{equation}
which, by taking the base-2 logarithm of both sides of each approximation, are equivalent to
\begin{equation} 
 ax - \log_2 (6/5) \approx 0, \qquad bx - \log_2 (5/4) \approx 0, \qquad (a+b)x - \log_2(3/2) \approx 0.
\label{apr3}
\end{equation}
To find the ``best'' $x$ that satisfies all three of the approximations in~\eqref{apr3} simultaneously,
we seek the absolute minimum of the function
\begin{equation} \label{min2param}
g(x) = \big(ax - \log_2 (6/5)\big)^2 +  \big(bx - \log_2 (5/4)\big)^2 + \big( (a+b)x - \log_2(3/2) \big)^2
\end{equation} over the positive real numbers.  Notice that $g(x)$ is the sum of the squared 
differences of the three approximations we seek to minimize; i.e., we will employ the 
principle of least squares to find our optimal $x$.

Since $g(x)$ is a quadratic function with a positive leading co\"efficient, it attains its absolute
minimum at $x = x_0$, where its derivative $g'(x_0) = 0$.  
Using elementary calculus, the $x$-value that minimizes $g(x)$ can be shown to be
\begin{equation} \label{CarlosUnit}
 x_0^{(a,b)} = \frac{1}{a^2 + b^2 + (a+b)^2} \log_2 \left[  \left( \frac 65 \right)^a \left( \frac 54 \right)^b 
    \left( \frac 32 \right)^{a+b} \right] .
\end{equation}

Let us refer to the scale with unit size given by $x_0^{(a,b)}$ in~\eqref{CarlosUnit} as 
the $(a,b)$-Carlos scale.  To express the unit size in cents, just multiply $x_0^{(a,b)}$ by $1200$.
The ratio of the frequency of one note to the next available lower note is 
\( 2^{x_0^{(a,b)}} . \)

Next, one needs to experiment with various values of $a$ and $b$, and check to see if 
all the
approximations are within an acceptable tolerance $v$, say, 
$v=$ $5$ cents. 
While $5$ cents may be considered an arbitrary cutoff and that 
additional scales of interest would emerge if we increased the tolerance to $8$ cents,
it seems that the conventional wisdom is that the trained human ear is sensitive to 
pitch differences of greater than or equal to $5$ cents; see, e.g.~\cite[p. 287, first full paragraph]{L87}, so this seems to be a reasonable, though perhaps arbitrary 
cutoff to employ for our present study.
Nonetheless, if we were to increase the limit to $5$, the fifths remain surprisingly stable and
always around $702$ cents. The thirds fluctuate a little more. Remember that
the distance from the perfect major third $(5/4)$ is $14$ cents for the tempered
major third $(2^{1/3})$ and $22$ cents for the Pythagorean major third $(81/64)$. And
the distance from the perfect minor third $(6/5)$ is about $16$ cents for the minor
third of the meantone temperament $(4/5^{3/4})$ or for the tempered minor third
and $22$ cents for the Pythagorean third $(32/27)$. In musical reality, the gaps
therefore seem much greater than in the theory of the article. If we go, for
example, from $v < 5$ to $v < 8$, we see that new scales appear, e.g., 
(11,13)-Carlos close to 41-ET, (15,19)-Carlos close to 48-ET.

The author programmed \emph{Mathematica} to perform such a search.
Recall that the just minor third, major third, and perfect fifth are 
$315.64$, $386.31$, and $701.96$ cents, respectively.  Thus, if
\begin{equation} \label{ErrorTolerance} 
\max \left\{  | 1200 a x_0^{(a,b)} - 315.64 |, | 1200 b x_0^{(a,b)} - 386.31 |, | 1200(a+b) x_0^{(a,b)} - 701.96| \right\} < 5, 
\end{equation}
then we have identified an \emph{$(a,b)$-Carlos scale} of interest.  

The minimal pair $(a,b)$ that satisfies~\eqref{ErrorTolerance} is $(4,5)$.  
The $(4,5)$-Carlos scale is the Carlos alpha scale, with
 \[ x_0^{(4,5)} = \frac{1}{122} \log_2 \left(  (6/5)^4 (5/4)^5 (3/2)^9 \right) = 
 \log_2 \left( \sqrt[122]{ \frac{7971615} {32768} }\right),\]
 and which is equivalent to a unit size of about 77.965 cents.   
 As is well known, an acceptably approximated octave does not exist in the alpha
 scale: it takes 15.3915 units to get a perfect octave, and 15 units falls about 30.5252 cents narrow.
 
 The $(5,6)$-Carlos scale is the Carlos beta scale, and the $(9,11)$-Carlos scale is the gamma scale.
 
 Another related scale noted by Carlos~\cite{C87} is the $\alpha'$ scale, which in our 
 notation is the
 $(8,10)$-Carlos scale.  It is simply the alpha scale with each unit
 cut in half (or equivalently, the number of units is doubled).  
 Here the unit is about $38.9825$ cents, half that of the alpha scale. 
 Of course, the thirds and the fifth are identical to that of the alpha scale; but this 
 ``doubled alpha'' scale has one advantage over the alpha scale: an octave requires about 30.783 units.  Thus the octave approximated by
 the 31st note is about 8.45734 cents wide, which is much better than the 15th note of the 
 alpha scale.  Nonetheless, 
 if $\gcd(a,b) > 1$, that $(a,b)$-Carlos scale could reasonably be considered a trivial variation
 of the $\big(a/\gcd(a,b), b/\gcd(a,b)\big)$-Carlos scale, and therefore not of genuine 
(mathematical)
 interest.\footnote{As usual, we denote the greatest common divisor of $a$ and
 $b$ by $\gcd(a,b)$.} 
 
Nonetheless, there is indeed \emph{musical} interest to consider
$(a,b)$-Carlos scales where $\gcd(a,b) > 1$.  For example, the $(8,10)$-Carlos scale,
called $\alpha'$ by Carlos, is close to 31-ET, the $(10,12)$-Carlos scale
($\beta'$) is close to 38-ET,  $(12,15)$-Carlos ($\alpha''$) is close to 48-ET, etc.
Further, one could ask the opposite question, namely for a given $N$-ET temperament,
which $a$ and $b$ give the closet $(a,b)$-Carlos approximation?  For example, the
H\"older temperament (53-ET) is well approximated by $(14,17)$-Carlos; see further 
comments on this below.

Carlos~\cite{C90} only considered potential scales with between $10$ and $40$ units
in her computer search.   Of course, $40$ is an arbitrary cutoff; Harry Partch, for example,
created a well known $43$-unit scale~\cite{P74}.   The 53-ET scale is also of 
particular importance; see, e.g.~\cite[p. 224 ff.]{B06}.
So let us consider possibilities beyond the $40$-unit per octave limit.
 
 Of particular interest is the $(14,17)$-Carlos scale.   Not only are the thirds and the fifth very close
 to just intonation, but serendipitously 53 units gives a note just 0.267329 cents wider than a 
 perfect octave!  Thus the $(14,17)$-Carlos scale is a very close approximation to 53-ET scale, but with slightly sweeter thirds and fifths, and an ever-so-slightly narrow octave. 
 
 Similarly of particular note for an excellent octave in addition to the thirds and the fifth is the
  $(17,21)$-Carlos scale,  as 65 units is 0.520 cents wider than an octave. 
   The $(19,23)$-Carlos scale in
 addition to excellent thirds and fifths also has a very pleasing octave: 72 units is 3.56 cents
wider than a perfect octave.  Please see the table below for additional examples of
 $(a,b)$-Carlos scales of interest, all with less than 81 units per octave.
The author does not intend to imply that there are no interesting examples with
 greater than 81 units per octave; the table just contains examples to suggest possibilities.
 For example, the $(25,31)$-Carlos scale is close to the
 sixteenth tone space (96-ET) often used in microtonal music. 
  \vskip 5mm
  
  \begin{tabular}{| c | r |  r | r |  r | r | c |}
  \hline 
  $(a,b)$ & unit size & m3 dev & M3 dev & P5 dev & units per  & notes\\
               &               & from just& from just &  from just &  octave    &   \\
  \hline
  $(4,5)$ & $77.965$ & $-3.780$ & $3.515$ & $-0.275$ & $15.39$ & alpha scale\\
  $(5,6)$ & $63.833$ & $3.525 $& $-3.312$ & $0.202$ & $18.80$ & beta scale \\
  $(9,11)$ & $35.099$ & $0.247$ & $-0.226$ & $0.011$ & $34.19$ & gamma scale \\
  $(13,16)$ & $24.203$ & $-1.000$ & $0.939$ & $-0.072$ &$49.58$ & \\
  $(14,17)$ & $22.647$ & $1.412$ & $-1.319$ & $0.083$ & $52.99$ & \\
  $(17,21)$ & $18.470$ & $-1.658$ & $1.550$ & $-0.118$ & $64.97$ & \\
  $(19,23)$ & $16.716$ & $-1.966$ & $-1.840$ & $0.116$ & $71.79$ & \\
  $(21,25)$ & $15.266$ & $4.945$ & $-4.661$ & $0.274$ & $78.61$ & \\
  $(21,26)$ & $14.932$ & $-2.064$ & $1.927$ & $-0.146$ & $80.36$ & \\
  \hline
  \end{tabular}
  \vskip 1cm
  For further discussion of the alpha, beta, and gamma scales, along with further thoughts
 on tuning and temperament by Carlos, see~\cite{C87}.
 
 \section{Further Generalizations with three or more parameters or with alternative intervals optimized}
 \subsection{Three or more independent parameters}
 Rather than requiring the perfect fifth to be exactly $a+b$ units, we could instead allow it to be $c$
 units where $a<b<c$.  Then the function to minimize is
 \begin{equation}\label{min3param}
 g(x) = \big(ax - \log_2 (6/5)\big)^2 +  \big(bx - \log_2 (5/4)\big)^2 + \big( cx - \log_2(3/2) \big)^2,
 \end{equation}
 and its minimum value is
 \begin{equation} \label{CarlosUnit3param}
 x_0^{(a,b,c)} = \frac{1}{a^2 + b^2 + c^2} \log_2 \left[  \left( \frac 65 \right)^a \left( \frac 54 \right)^b 
    \left( \frac 32 \right)^{c} \right] .
\end{equation}
We will name a scale built from this unit $x_0^{(a,b,c)}$ an $(a,b,c)$-Carlos scale. 
If it so happens that $c = a+b$, then this is just the $(a,b)$-Carlos scale above. 
But notice that here no attempt is made to control the size of the minor seventh as $a+c$
units nor the major seventh as $b+c$ units.

 An example of such a scale can be found in the Wikipedia article on the ``Delta scale"
 \cite{WDelta}: what is
 described there is the $(23,28,50)$-Carlos scale in our nomenclature.  Notice that 
 $23 + 28 \neq 50$, and therefore is not the same as the $(23,28)$-Carlos scale defined above.

 It is worth 
 acknowledging here the shortcomings of scales where the approximation to $6/5$ plus
 the approximation to $5/4$ does not exactly equal the interval approximating $3/2$; and
 that this phenomenon is called \emph{inconsistency}, and can make traditional
 composing techniques with chords awkward.
 
 \subsection{Alternative intervals optimized}
 Carlos decided to optimize the intonation of the major and minor triads to the possible detriment of
 the octave, thus giving equal unit scales with very near just minor and minor thirds and perfect 
 fifths.  But one could alternatively choose to optimize other intervals.  For example, if we 
 allocate $a$ units for the perfect fourth and $b$ units for the perfect fifth, and therefore $a+b$
 units for the perfect octave, the analog of Eq.~\eqref{min2param} becomes
 \begin{equation}\label{min2paramAlt}
 g(x) = \big(ax - \log_2 (4/3)\big)^2 +  \big(bx - \log_2 (3/2)\big)^2 + \big( (a+b)x - 1 \big)^2.
 \end{equation}
 The unique minimizer of~\eqref{min2paramAlt} is
 \begin{align} \label{CarlosUnitP4P5}
 x_{ \{ \mathrm{P4,P5}\} }^{(a,b)}
 &= \frac{1}{a^2 + b^2 + (a+b)^2} \log_2 \left[  \left( \frac 43 \right)^a \left( \frac 32 \right)^b 
    2^{a+b} \right] \\
  &=    \frac{3a + (b-a)\log_2 3}{a^2 + b^2 + (a+b)^2} , \notag
\end{align}
where we introduce the notation  $x_{ \{ {I_a,I_b}\} }^{(a,b)}$ to denote the
unit size of the Carlos-type scale where the just interval $I_k$ is approximated by $k$ units, for
$k = a$ or $b$, and the interval $I_a \oplus I_b$ is approximated
by $a+b$ units.  Addition of intervals is defined in the intuitive way, e.g. 
$\mathrm{m3\oplus M3 = P5}$,
$\mathrm{P4 \oplus P5 = P8}$, 
$\mathrm{M2 \oplus m3 = P4}$, etc.   In this way, the previously introduced
$ x_{0}^{(a,b)}$ is the same as  $x_{ \{ \mathrm{m3,M3}\} }^{(a,b)}$ in this more general
notation.
 
 Note that unit size  $x_{ \{ \mathrm{P4,P5}\} }^{(5,7)} 
 = (15 + 2\log_2 3)/218 \approx 0.0833483 \approx 0.08333333\dots = 1/12$ gives a very close approximation to 12-ET.  
One could argue, however, that the preceding example is an abuse of
Carlos' original idea, as it  could be considered 
wasteful to optimize for both the perfect fourth and fifth, as each is the
inversion of the other.  Indeed, the most simple solution $(a,b) = (5,7)$ does not really get us
anything new, but rather a slightly out of tune version of 12-ET(!) 
 
 Letting $J_{I}$ denote the just ratio for interval $I$, so that, e.g. $J_{\mathrm{P5}} = 3/2$,
 $J_{\mathrm{M3}} = 5/4$, etc., we can produce scales with optimized approximations for
 intervals $I_a$, $I_b$, and $I_a \oplus I_b$ using unit
\begin{equation} \label{CarlosUnitGen}
 x_{ \{ I_a, I_b\} }^{(a,b)} = \frac{1}{a^2 + b^2 + (a+b)^2} 
    \log_2 \left[  (J_{I_a})^a  (J_{I_b})^b  (J_{I_a \oplus I_b})^{a+b} \right] .
\end{equation}
To be careful, we may choose to exclude possibilities where 
intervals $I_a$ and $I_b$ are inversions of
each other.

\subsection{Example: a Carlos-type pentatonic optimized scale}
We will provide one more example to demonstrate the flexibility of the proposed optimization scheme
even with only two independent parameters $a < b$.
Consider a pentatonic scale given by the consecutive intervals M2, M2, m3, M2. 
(Relative to the tonic, the intervals are M2, M3, P5, M6.)
We wish to set up a Carlos-type scale with unit size determined so that the intervals
in the pentatonic scale (rather than the major and minor triads) are as close to just
as possible.
 Letting the major second be
$a$ units apart and the minor third $b$ units apart, we would then want the major third to be
$2a$ units, the perfect fifth to be $2a+b$ units.   
(In light of the earlier remark, we exclude the condition that the major sixth ought to be 
approximated by $3a+b$ units, as the major sixth is
the inversion of the already included minor third.)
Thus, following from the discussion earlier, we want find the unique minimizer of the function
\begin{equation} \label{min2paramPent}
g(x) = 
\big(ax - \log_2  (9/8) \big)^2 +  \big(2ax - \log_2 (5/4)\big)^2 +
 \big( bx - \log_2(6/5) \big)^2
 +\big( (2a+b)x - \log_2 (3/2) \big)^2
\end{equation} over the positive real numbers, which is unit size
\begin{multline}
 \frac{1}{a^2 + (2a)^2 + b^2 + (2a+b)^2} \log_2 \left[ 
     \left( \frac 98 \right)^a
     \left( \frac 54 \right)^{2a}
     \left( \frac 65 \right)^b
     \left( \frac 32 \right)^{2a+b}
    \right]\\ = 
   \frac{1}{9a^2 + 4ab + 2b^2} \log_2 \left(  \frac{5^{2a-b} 9^{2a+b}}{512^a}  \right) . 
\end{multline} 
For $(a,b)=(7,11)$, we obtain a scale where each of the major second,
minor third, major third, and perfect fifth differs from its just counterpart by less than
7.25 cents.  The 43-unit octave is 8.71617 cents wide.

Another notable Carlos-type pentatonic-optimized scale occurs at $(a,b) = (17,27)$, where
each of the major second, minor third, major third, and perfect fifth differs from its
just counterpart by less than 8 cents, and the 104-unit octave is only 1.3501 cents
narrow.

\section{Conclusion}
From the examples above, interested readers should be able to construct many 
additional variants that will fall under the general umbrella of Carlos-type scales.

\section*{Acknowledgments}
The author wishes to thank Robert Schneider and Eden Sills for their assistance and encouragement, 
and for valuable suggestions.  The author also thanks the
Mathematics and Music Lab at Michigan Technological University for invaluable 
consultation work.   Finally, the author thanks the anonymous referees for catching a
typographical error in the original submission, and for making many useful suggestions
that improved the paper considerably.

\section*{Disclosure statement}
The author reports that there are no competing interests to declare.

\end{document}